\documentclass[12pt]{article}
\usepackage{amssymb}
\usepackage{amsmath}

\newcommand{\lap}{\mbox{$\bigtriangleup$}}
\newcommand{\grad}{\mbox{$\bigtriangledown$}}
\newcommand{\ra}{{\mbox{$\rightarrow$}}}
\newcommand{\be}{\begin{equation}}
\newcommand{\ee}{\end{equation}}

\newtheorem{mthm}{Theorem}

\newtheorem{thm}{Theorem}[section]

\begin{document}

\title{Maximum principles for a fully nonlinear fractional order equation and symmetry of solutions}

\author{Wenxiong Chen \thanks{Partially supported by the Simons Foundation Collaboration Grant for Mathematicians 245486.} \quad  Congming Li \thanks{ Corresponding author, Partially supported by NSF DMS-1405175 and NSFC-11271166.} \quad Guanfeng Li  }

\date{\today}
\maketitle
\begin{abstract}
In this paper, we consider equations involving fully nonlinear nonlocal operators
$$F_{\alpha}(u(x)) \equiv C_{n,\alpha} PV \int_{\mathbb{R}^n} \frac{G(u(x)-u(z))}{|x-z|^{n+\alpha}} dz= f(x,u).$$

We prove a {\em maximum principle} and obtain key ingredients
for carrying on the method of moving planes, such as {\em narrow region principle} and {\em decay at infinity}. Then we establish radial symmetry and monotonicity for positive solutions to Dirichlet problems associated to such fully nonlinear fractional order equations in a unit ball and in the whole space, as well as non-existence of solutions on a half space. We believe that the methods develop here
can be applied to a variety of problems involving fully nonlinear nonlocal operators.

We also investigate the limit of this operator as $\alpha \ra 2$ and show that
$$F_{\alpha}(u(x)) \ra a(-\lap u(x)) + b |\grad u(x)|^2 .$$

\end{abstract}
\bigskip

\bigskip

\section{Introduction}

In this paper, we consider nonlinear equations involving fully nonlinear nonlocal operators
\begin{equation}
F_{\alpha}(u) = f(x,u)
\label{Feq}
\end{equation}
with
\begin{eqnarray}
F_{\alpha}(u(x)) &=& C_{n,\alpha} \, \lim_{\epsilon \ra 0} \int_{\mathbb{R}^n\setminus B_{\epsilon}(x)} \frac{G(u(x)-u(z))}{|x-z|^{n+\alpha}} dz \nonumber \\
&=& C_{n,\alpha} PV \int_{\mathbb{R}^n} \frac{G(u(x)-u(z))}{|x-z|^{n+\alpha}} dz,
\label{F}
\end{eqnarray}
where PV stands for the Cauchy principle value. This kinds of operator was introduced by Caffarelli and Silvestre in \cite{CS1}).

In order the integral to make sense, we require that
$$u \in C^{1,1}_{loc} \cap L_{\alpha} $$
with
$$ L_{\alpha} = \{ u \in L^1_{loc} \mid \int_{\mathbb{R}^n} \frac{|u(x)|}{1+|x|^{n+\alpha}} d x < \infty \};$$
$G$ being at least local Lipschitz continuous, and $G(0)=0$.

In the special case when $G(\cdot)$ is an identity map, $F_{\alpha}$ becomes the fractional Laplacian $(-\lap)^{\alpha/2}$. The nonlocal nature of these operators make them difficult to study. To circumvent this, Caffarelli and Silvestre \cite{CS} introduced the {\em extension method } which turns the nonlocal problem involving the fractional Laplacian into a local one in higher dimensions. This method has been applied successfully to study equations involving the fractional Laplacian, and a series of fruitful results have been obtained (see \cite{BCPS} \cite{CZ}    and the references therein). One can also use the {\em integral equations method}, such as {\em the method of moving planes in integral forms} and {\em regularity lifting} to investigate equations involving fractional Laplacian by first showing that they are equivalent to the corresponding integral equations. \cite{CLO} \cite{CLO1} \cite{CFY}.

However, for the fully nonlinear nonlocal operator $F_{\alpha}(\cdot)$, so far as we know, there has neither  been any corresponding {\em extension methods} nor equivalent integral equations that one can work at. This is probably the reason that very few results has been obtained for such fully nonlinear nonlocal operators. Hence it is essential to develop methods that can deal directly on these kinds of nonlocal operators, which is the main objective of our paper.

We first prove

\begin{mthm} (The Simple Maximum Principle)

Let $\Omega$ be a bounded domain in $\mathbb{R}^n$. Assume that $u \in C^{1,1}_{loc} (\Omega) \cap L_{\alpha}$, be lower semi-continuous on $\bar{\Omega}$, and satisfies
\begin{equation}
\left\{\begin{array}{ll}
F_{\alpha}(u(x)) \geq 0, & x \in \Omega,\\
u(x) \geq 0, & x \in \Omega^c.
\end{array}
\right.
\end{equation}
Suppose
$$G(0)= 0, \;\; G \in C^1(\mathbb{R}), \; \mbox{ and } G'(t) \geq c_0 > 0, \forall t \in \mathbb{R}.$$

Then
\begin{equation}
u(x) \geq 0, \;\;\;\; x \in \Omega.
\end{equation}

The same conclusion holds for unbounded domains $\Omega$ if we further assume that
$$\liminf_{|x| \ra \infty} u(x) \geq 0.$$
\label{mthm1}
\end{mthm}

Then we establish maximum principles for anti-symmetric functions which play important roles in carrying out the method of moving planes. To explain the assumptions in these principles, we recall some basics in this method.

Take the whole space $\mathbb{R}^n$ as an example. Let
$$T_{\lambda} =\{x \in \mathbb{R}^{n}|\; x_1=\lambda, \mbox{ for some } \lambda\in \mathbb{R}\}$$
be the moving planes,
$$\Sigma_{\lambda} =\{x \in \mathbb{R}^{n} | \, x_1<\lambda\}$$
be the region to the left of the plane, and
$$ x^{\lambda} =(2\lambda-x_1, x_2, ..., x_n)$$
be the reflection of $x$ about the plane $T_{\lambda}$.

Assume that $u$ is a solution of pseudo differential equation (\ref{Feq}). To compare the values of $u(x)$ with $u(x^{\lambda})$, we denote
$$w_{\lambda} (x) = u(x^{\lambda}) - u(x) .$$

The first step is to show that for $\lambda$ sufficiently negative, we have
\begin{equation}
w_{\lambda}(x) \geq 0 , \;\; x \in \Sigma_{\lambda} .
\label{w}
\end{equation}
This provides a starting point to move the plane. Then in the second step, we move the plane to the right as long as inequality (\ref{w}) holds to its limiting position to show that $u$ is symmetric
about the limiting plane. Usually, a maximum principle is used to prove (\ref{w}).
From (\ref{Feq}), we have
$$ F(u_{\lambda}) - F(u) = f(x, u_{\lambda}) - f(x,u) = - c_{\lambda}(x) w_{\lambda}(x) ,$$
for some $c_{\lambda}(x)$ depending on $\lambda$ and $u(x)$.
Under appropriate assumptions on $f$ and $u$, or after making a Kelvin transform, this $c_{\lambda}(x)$ may have a certain rate of decay near infinity;
and it is easy to see that $w_{\lambda}$ is an anti-symmetric function:
$$w_{\lambda}(x) = - w_{\lambda}(x^{\lambda}).$$

In order to obtain (\ref{w}) in step 1, we establish a {\em decay at infinity principle}.
For simplicity of notation, in the following, we denote $w_{\lambda}$ by $w$ and $\Sigma_{\lambda}$ by $\Sigma$.

\begin{mthm} ( Decay at Infinity.)

Let $\Omega$ be an unbounded region in $\Sigma$.
Assume $w \in L_{\alpha} \cap C_{loc}^{1,1}(\Omega)$ is a solution of
$$
\left\{\begin{array}{ll}
F(u_{\lambda}(x))-F(u(x)) +c(x)w(x)\geq0  &\mbox{ in } \Omega,\\
w(x) \geq0&\mbox{ in } \Sigma \backslash\Omega,\\
w(x^{\lambda})=-w(x)  &\mbox{ in } \Sigma,
\end{array}
\right.
$$
with
\begin{equation}
\underset{|x|\rightarrow \infty}{\underline{\lim}}|x|^{\alpha}c(x)\geq0,
\label{cx}
\end{equation}
then there exists a constant $R_0>0$ ( depending on $c(x)$, but independent of $w$ ), such that if
$$
w(x^0)=\underset{\Omega}{\min}\;w(x)<0,
$$
 then
 $$
|x^0|\leq R_0.
$$
\label{mthm2}
\end{mthm}

Note that the condition (\ref{cx}) is satisfied if
$$
c(x) = o(\frac{1}{|x|^{\alpha}}), \;\; \mbox{ for $|x|$ sufficiently large}.
$$

In the second step, let
$$\lambda_o = \{ \lambda \mid w_{\mu}(x) \geq 0, \; x \in \Sigma_{\mu}, \, \mu \leq \lambda\}$$
be the upper limit of such $\lambda$ that (\ref{w}) holds.
To show that $u$ is symmetric about the limiting plane $T_{\lambda_o}$, or
\begin{equation}
 w_{\lambda_o}(x) \equiv 0 , \;\; x \in \Sigma_{\lambda_o};
 \label{w0}
 \end{equation}
we usually use a contradiction argument: if (\ref{w0}) does not hold, then we can move the
plane a little bit forward, and still have (\ref{w}) for some $\lambda > \lambda_o$. The region
between $T_{\lambda_o}$ and $T_{\lambda}$ is a narrow region, and the {\em maximum principle} holds in a narrow region provided $c(x)$ is not ``too negative'', as you will see below.

\begin{mthm} ( Narrow Region Principle.)

Let $\Omega$ be a bounded narrow region in $\Sigma$, such that it is contained in
 $$\{x | \; \lambda-\delta<x_1<\lambda \, \}$$
 with small $\delta$. Suppose  that $w\in L_\alpha\cap C_{loc}^{1,1}(\Omega)$ and is lower semi-continuous on $\bar{\Omega}$. If
 $c(x)$ is bounded from below in $\Omega$  and
$$
\left\{\begin{array}{ll}
F(u_{\lambda}(x))-F(u(x)) +c(x)w(x)\geq0  &\mbox{ in } \Omega,\\
w(x) \geq0&\mbox{ in }  \Sigma \backslash\Omega,\\
w(x^{\lambda})=-w(x)  &\mbox{ in } \Sigma,
\end{array}
\right.
$$

then for sufficiently small $\delta$, we have
$$
w(x) \geq0 \mbox{ in } \Omega.
$$

Furthermore, if $w = 0$ at some point in $\Omega$, then
$$ w(x) = 0 \; \mbox{ almost everywhere in }  \mathbb{R}^n. $$

These conclusions hold for unbounded region $\Omega$ if we further assume that
$$\underset{|x| \ra \infty}{\underline{\lim}} w(x) \geq 0 .$$
\label{mthm3}
\end{mthm}

We will use several examples to illustrate how the key ingredients in the above can be used in the {\em method of moving planes} to establish symmetry and monotonicity of positive solutions.

We first consider
\be
\left\{ \begin{array}{ll}
F_{\alpha}(u(x))  = f(u(x)) , & x \in B_1(0), \\
u(x) = 0 , & x \not{\in} B_1(0).
\end{array}
\right.
\label{mfu1}
\ee
We prove
\begin{mthm}
Assume that $u \in L_{\alpha}\cap C_{loc}^{1,1}(B_1(0))$ is a positive solution of (\ref{mfu1}) with $f(\cdot)$ being Lipschitz continuous. Then $u$ must be radially symmetric and monotone decreasing about the origin.
\label{mthm4}
\end{mthm}

Then we study
\be
F_{\alpha}(u(x)) = g(u(x)), \;\; x \in \mathbb{R}^n .
\label{meqws}
\ee

\begin{mthm}
Assume that $u \in C^{1,1}_{loc} \cap L_{\alpha}$ is a positive solution of (\ref{meqws}).
Suppose, for some $\gamma >0$,
$$
u(x)=o(\frac{1}{|x|^{\gamma}}), \mbox{ as } |x|\ra \infty,
$$
and
$$
g'(s) \leq s^q , \;\; \mbox{ with } q \gamma \geq \alpha .
$$
Then $u$ must be radially symmetric about some point in $\mathbb{R}^n$.
\label{mthm5}
\end{mthm}

Finnaly, we investigate a Dirichlet problem on an upper half space
$$\mathbb{R}_+^n = \{ x = (x_1, \cdots, x_n) \mid x_n > 0 \}.$$

Consider
\begin{equation}
\left\{\begin{array}{ll}
F_{\alpha}(u(x))= h(u(x)), & x \in \mathbb{R}_+^n, \\
u(x) \equiv 0 , & x \not{\in} \mathbb{R}_+^n .
\end{array}
\right.
\label{me1h1}
\end{equation}

\begin{mthm}
Assume that $u \in L_{\alpha} \cap C^{1,1}_{loc}$ is a nonnegative solution of
problem (\ref{me1h1}). Suppose
$$ \lim_{|x| \ra \infty} u(x) = 0,
$$
$$
h(s) \mbox{ is Lipschitz continuous in the range of } u, \mbox{ and } h(0)=0.
$$
Then $u \equiv 0$.
\label{mthm6}
\end{mthm}

Finally, we investigate the limit  of $F_{\alpha}(u(x))$ as $\alpha \ra 2$ for each fixed $x$ and
discover an interesting phenomenon.

\begin{mthm}
Assume that $u \in C^{1,1}_{loc} \cap L_{\alpha}$, $G(\cdot)$ is second order differentiable, and $G(0) = 0$.
Then
$$ \lim_{\alpha \ra 2} F_{\alpha}(u(x)) = a (-\lap u)(x) + b |\grad u(x)|^2 ,
$$
where $a$ and $b$ are constant multiple of $G'(0)$ and $G''(0)$ respectively.
\label{mthm7}
\end{mthm}

In Section 2, we establish various maximum principles and prove Theorem \ref{mthm1}, \ref{mthm2}, and \ref{mthm3}. In Section 3, we use the key ingredients obtained in Section 2 to derive symmetry and nonexistence of solutions and prove Theorem \ref{mthm4}, \ref{mthm5}, and \ref{mthm6}.
In Section 4, we derive Theorem \ref{mthm7}.

For more articles concerning the method of moving planes for nonlocal equations and for integral equations, please see \cite{FL} \cite{Ha} \cite{HLZ} \cite{HWY} \cite{Lei} \cite{LLM} \cite{LZ} \cite{LZ1} \cite{LZ2} \cite{MC} \cite{MZ} and the references therein.

\section{Various Maximum Principles}
Let
\begin{equation}
F_{\alpha}(u(x)) = C_{n,\alpha} \, PV \int_{\mathbb{R}^n} \frac{G(u(x)-u(z))}{|x-z|^{n+\alpha}} dz.
\label{F1}
\end{equation}

Throughout this and the next section, we assume that $G \in C^1(\mathbb{R})$,
\begin{equation}
G(0) =0, \; \mbox{ and } \; G'(t) \geq c_0 > 0 \;\; \forall \, t \in \mathbb{R}.
\label{2.1}
\end{equation}

\begin{thm} (The Simple Maximum Principle)

Let $\Omega$ be a bounded domain in $\mathbb{R}^n$. Assume that $u \in C^{1,1}_{loc} (\Omega) \cap L_{\alpha}$, be lower semi-continuous on $\bar{\Omega}$, and satisfies
\begin{equation}
\left\{\begin{array}{ll}
F_{\alpha}(u(x)) \geq 0, & x \in \Omega,\\
u(x) \geq 0, & x \in \Omega^c.
\end{array}
\right.
\label{2.2}
\end{equation}
Then
\begin{equation}
u(x) \geq 0, \;\;\;\; x \in \Omega.
\label{2.3}
\end{equation}

The same conclusion holds for unbounded domains $\Omega$ if we further assume that
$$\liminf_{|x| \ra \infty} u(x) \geq 0.$$
\label{thm2.1}
\end{thm}

{\bf Proof.} Suppose (\ref{2.3}) is violated, then since $u$ is lower semi-continuous on $\bar{\Omega}$, there exists $x^o$ in $\Omega$ such that
$$u(x^o) = \min_{\Omega} u < 0.$$
It follows from (\ref{2.1}) that
\begin{eqnarray*}
\int_{\mathbb{R}^n} \frac{G(u(x^o)-u(z))}{|x^o-z|^{n+\alpha}} dz
&=& \int_{\mathbb{R}^n} \frac{G'(\psi(z))[u(x^o)-u(z)]}{|x^o-z|^{n+\alpha}} dz\\
&=& \leq c_0 \int_{\mathbb{R}^n} \frac{[u(x^o)-u(z)]}{|x^o-z|^{n+\alpha}} dz\\
&<& 0.
\end{eqnarray*}
This contradicts (\ref{2.2}) and hence proves the theorem.

In the following, we will continue to use the notation introduced in the previous section. Let

$$T_{\lambda} =\{x \in \mathbb{R}^{n}|\; x_1=\lambda, \mbox{ for some } \lambda\in \mathbb{R}\}$$
be the moving planes,
$$\Sigma_{\lambda} =\{x \in \mathbb{R}^{n} | \, x_1<\lambda\}$$
be the region to the left of the plane, and
$$ x^{\lambda} =(2\lambda-x_1, x_2, ..., x_n)$$
be the reflection of $x$ about the plane $T_{\lambda}$.
$$w_{\lambda} (x) = u(x^{\lambda}) - u(x) .$$
For simplicity of notation, we denote $w_{\lambda}$ by $w$ and $\Sigma_{\lambda}$ by $\Sigma$.

\begin{thm}(Maximum Principle for Anti-symmetric Functions)
Let $\Omega$ be a bounded domain in $\Sigma$.
Assume that $w\in L_\alpha\cap C_{loc}^{1,1}(\Omega)$ and is lower semi-continuous on $\bar{\Omega}$. If
$$
\left\{\begin{array}{ll}
F(u_{\lambda}(x)) - F(u(x))  \geq 0  &\mbox{ in } \Omega,\\
w(x) \geq0&\mbox{ in }  \Sigma \backslash\Omega,\\
w(x^{\lambda})=-w(x)  &\mbox{ in } \Sigma,
\end{array}
\right.
$$
then
$$
w(x) \geq0 \mbox{ in } \Omega.
$$

Furthermore, if $w = 0$ at some point in $\Omega$, then
$$ w(x) = 0 \; \mbox{ almost everywhere in }  \mathbb{R}^n. $$

These conclusions hold for unbounded region $\Omega$ if we further assume that
$$\underset{|x| \ra \infty}{\underline{\lim}} w(x) \geq 0 .$$

\label{thm2.2}
\end{thm}

{\bf Proof.} Suppose otherwise, then there exists a point $x$ in $\Omega$, such that
\begin{equation}
w(x) = \min_{\Omega} w < 0 .
\label{w<0}
\end{equation}

\begin{eqnarray}
F_{\alpha}(u_{\lambda}(x)) - F_{\alpha}(u(x)) &=& C_{n, \alpha} PV \int_{\mathbb{R}^n} \frac{G(u_{\lambda}(x) -u_{\lambda}(y)) - G(u(x)-u(y))}{|x-y|^{n+\alpha}} d y \nonumber \\
&=& C_{n, \alpha} PV \int_{\Sigma} \frac{G(u_{\lambda}(x) -u_{\lambda}(y)) - G(u(x)-u(y))}{|x-y|^{n+\alpha}} d y\nonumber \\
 &+& C_{n, \alpha} PV \int_{\Sigma} \frac{G(u_{\lambda}(x) -u(y)) - G(u(x)-u_{\lambda}(y))}{|x-y^{\lambda}|^{n+\alpha}} d y \nonumber\\
&=& C_{n, \alpha} PV \left\{ \int_{\Sigma_-} [ \cdots ] d y +  \int_{\Sigma_+} [ \cdots ] d y \right\} \nonumber \\
&=& C_{n, \alpha} PV \left\{ I_1 + I_2 \right\}.
\label{2.5}
\end{eqnarray}

Here $$\Sigma_- = \{ y \in \Sigma \mid w(x) +w(y) \leq 0 \} $$
and $$\Sigma_+ = \{ y \in \Sigma \mid w(x) +w(y) > 0 \} .$$

To estimate $I_1$, we use (\ref{2.1}) and the fact
$$ \frac{1}{|x-y|} > \frac{1}{|x-y^{\lambda}|} , \;\; \forall \, x, y \in \Sigma $$
to arrive at
\begin{eqnarray}
I_1 &\leq& \int_{\Sigma_-} \left\{ \frac{G'(\cdot)(w(x)-w(y))}{|x-y|^{n+\alpha}} + \frac{G'(\cdot)(w(x)+w(y))}{|x-y^{\lambda}|^{n+\alpha}} \right\} d y \nonumber \\
&\leq& \int_{\Sigma_-} \left\{ \frac{c_0 (w(x)-w(y))}{|x-y|^{n+\alpha}} + \frac{c_0 (w(x)+w(y))}{|x-y^{\lambda}|^{n+\alpha}} \right\} d y \nonumber \\
&\leq& \int_{\Sigma_-} \left\{ \frac{c_0 (w(x)-w(y))}{|x-y^{\lambda}|^{n+\alpha}} + \frac{c_0 (w(x)+w(y))}{|x-y^{\lambda}|^{n+\alpha}} \right\} d y \nonumber \\
&\leq & 2 c_0 w(x) \int_{\Sigma_-} \frac{1}{|x-y^{\lambda}|^{n+\alpha}} d y .
\label{2.6}
\end{eqnarray}

To estimate $I_2$, we first notice that on $\Sigma_+$,  by (\ref{2.1}), we have
$$G(u_{\lambda}(x) - u(y)) - G(u(x)-u_{\lambda}(y))= G'(\cdot)(w(x)+w(y)) > 0, $$
and hence
\begin{eqnarray}
& & I_2 \nonumber \\
&\leq& \int_{\Sigma_+} \frac{G(u_{\lambda}(x) -u_{\lambda}(y)) - G(u(x)-u(y)) +G(u_{\lambda}(x) -u(y)) - G(u(x)-u_{\lambda}(y))}{|x-y|^{n+\alpha}} d y\nonumber \\
&=& \int_{\Sigma_+} \frac{[G(u_{\lambda}(x) -u_{\lambda}(y)) - G(u(x)-u_{\lambda}(y))] +[G(u_{\lambda}(x) -u(y)) - G(u(x)-u(y))] }{|x-y|^{n+\alpha}} d y \nonumber \\
&=& \int_{\Sigma_+} \frac{G'(\cdot) w(x) + G'(\cdot) w(x) }{|x-y|^{n+\alpha}} d y \nonumber \\
&=& 2 c_0 w(x) \int_{\Sigma_+} \frac{1}{|x-y|^{n+\alpha}} d y .
\label{2.7}
\end{eqnarray}

Combining (\ref{2.5}), (\ref{2.6}), and (\ref{2.7}), we derive
\begin{eqnarray}
F_{\alpha}(u_{\lambda}(x)) - F_{\alpha}(u(x)) &\leq& 2 C_{n, \alpha}  c_0 w(x) \left\{ \int_{\Sigma_-} \frac{1}{|x-y^{\lambda}|^{n+\alpha}} d y + \int_{\Sigma_+} \frac{1}{|x-y|^{n+\alpha}} d y \right\} \nonumber \\
&\leq& 2 C_{n, \alpha}  c_0 w(x) \left\{ \int_{\Sigma_-} \frac{1}{|x-y^{\lambda}|^{n+\alpha}} d y + \int_{\Sigma_+} \frac{1}{|x-y|^{n+\alpha}} d y \right\} \nonumber \\
&=&  2 C_{n, \alpha}  c_0 w(x) \int_{\Sigma} \frac{1}{|x-y^{\lambda}|^{n+\alpha}} d y  \label{ki}\\
&<& 0. \nonumber
\end{eqnarray}
This is a contradiction with the equation and hence we must have
\begin{equation}
w_{\lambda}(x) \geq 0, \;\; \forall \, x \in \Omega.
\label{2.1b}
\end{equation}

{\em Note inequality (\ref{ki}) we just derived holds at a negative minimum $x$ of $w$, and it is a key ingredient in obtaining various {\em maximum principles},
which will be used several times later.}

Suppose $w=0$ at some point in $\Omega$, say $w(x^o)=0$, then we must have
\begin{equation}
w(x) = 0 \; \mbox{ almost everywhere in }  \mathbb{R}^n.
\label{2.2b}
\end{equation}

To this end, we re-estimate $I_1$ and $I_2$ at $x=x^o$. From the previous argument, we have

\begin{eqnarray}
I_1  &\leq& \int_{\Sigma_-} \left\{ \frac{c_0 (w(x^o)-w(y))}{|x^o-y|^{n+\alpha}} + \frac{c_0 (w(x^o)+w(y))}{|x^o-y^{\lambda}|^{n+\alpha}} \right\} d y \nonumber \\
&=& c_0 \int_{\Sigma_-} \left( \frac{1}{|x^o-y^{\lambda}|^{n+\alpha}} - \frac{1}{|x^o-y|^{n+\alpha}} \right) w(y) dy.
\label{2.4b}
\end{eqnarray}

\begin{eqnarray}
&&I_2 = \int_{\Sigma_+} \left( \frac{1}{|x^o-y^{\lambda}|^{n+\alpha}} - \frac{1}{|x^o-y|^{n+\alpha}} \right) \left[ G(u_{\lambda}(x^o)-u(y)) - G(u(x)-u_{\lambda}(y)) \right] d y\nonumber \\
&+& \int_{\Sigma_+} \frac{G(u_{\lambda}(x^o) -u_{\lambda}(y)) - G(u(x^o)-u(y)) +G(u_{\lambda}(x^o) -u(y)) - G(u(x^o)-u_{\lambda}(y))}{|x^o-y|^{n+\alpha}} d y\nonumber \\
&\leq& \int_{\Sigma_+} \left( \frac{1}{|x^o-y^{\lambda}|^{n+\alpha}} - \frac{1}{|x^o-y|^{n+\alpha}} \right) G'(\cdot) [ w(x^o)+w(y) ] d y + 2 c_0 w(x^o) \int_{\Sigma_+} \frac{1}{|x-y|^{n+\alpha}} d y . \nonumber\\
&\leq& c_0 \int_{\Sigma_+} \left( \frac{1}{|x^o-y^{\lambda}|^{n+\alpha}} - \frac{1}{|x^o-y|^{n+\alpha}} \right) w(y) d y.
\label{2.3b}
\end{eqnarray}

Combining (\ref{2.5}), (\ref{2.4b}), and (\ref{2.3b}), we derive
$$F_{\alpha}(u_{\lambda}(x^o)) - F_{\alpha}(u(x^o)) \leq C_{n, \alpha} c_0 PV \int_{\Sigma} \left( \frac{1}{|x^o-y^{\lambda}|^{n+\alpha}} - \frac{1}{|x^o-y|^{n+\alpha}} \right) w(y) d y.$$

Notice that the difference in the above parentheses is strictly negative while $w(y)$ is nonnegative, if
$w(y)$ is not identically zero in $\Sigma$, this would contradicts
$$ F_{\alpha}(u_{\lambda}(x^o)) - F_{\alpha}(u(x^o)) \geq 0 .$$
Therefore, we must have
$$w(y) \equiv 0, \;\; \forall \, y \in \Sigma, $$
and consequently, by symmetry, we arrive at (\ref{2.2b}).

If $\Omega$ is an unbounded region, then under the condition
$$\underset{|x| \ra \infty}{\underline{\lim}} w(x) \geq 0 ,$$
any negative minimum of $w$ in $\Omega$ is attained at some point $x^o \in \Omega$, and similar to the
previous arguments, we can deduce all the same conclusions.

This completes the proof of the theorem.

In many cases, the inequality
$$ F(u_{\lambda}(x)) - F(u(x)) \geq 0 $$
may not be satisfied as required in the previous theorem. However one can derive that
$$ F(u_{\lambda}(x)) - F(u(x))  + c(x) w (x) \geq 0 $$
for some function $c(x)$ depending on $u(x)$ and $\lambda$. If $c(x)$ is nonnegative, it is easy to see that the {\em maximum principle} is still valid; however this is usually not the case in practice.
Fortunately, in the process of moving planes, each time we only need to move $T_{\lambda}$ a little bit forward, hence the increment of $\Sigma_{\lambda}$ is a narrow region, and a {\em maximum principle} is easier to hold in a narrow region provided $c(x)$ is not ``too negative'', as you will see below.

\begin{thm} ( Narrow Region Principle.)

Let $\Omega$ be a bounded narrow region in $\Sigma$, such that it is contained in
 $$\{x | \; \lambda-\delta<x_1<\lambda \, \}$$
 with small $\delta$. Suppose  that $w\in L_\alpha\cap C_{loc}^{1,1}(\Omega)$ and is lower semi-continuous on $\bar{\Omega}$. If
 $c(x)$ is bounded from below in $\Omega$  and
$$
\left\{\begin{array}{ll}
F(u_{\lambda}(x)) - F(u(x)) +c(x)w(x)\geq0  &\mbox{ in } \Omega,\\
w(x) \geq0&\mbox{ in }  \Sigma \backslash\Omega,\\
w(x^{\lambda})=-w(x)  &\mbox{ in } \Sigma,
\end{array}
\right.
$$

then for sufficiently small $\delta$, we have
$$
w(x) \geq0 \mbox{ in } \Omega.
$$

Furthermore, if $w = 0$ at some point in $\Omega$, then
$$ w(x) = 0 \; \mbox{ almost everywhere in } \mathbb{R}^n. $$

These conclusions hold for unbounded region $\Omega$ if we further assume that
$$\underset{|x| \ra \infty}{\underline{\lim}} w(x) \geq 0 .$$
\label{thm2.3}
\end{thm}

{\bf Proof.}

Suppose in the contrary, there exists an $x^0 \in \Omega$, such that
$$w(x^0) = \min_{\Omega} w < 0. $$
Then by the {\em key inequality} (\ref{ki}), we deduce
\begin{eqnarray}
& & F(u_{\lambda}(x^0)) - F(u(x^0)) + c(x^0) w(x^0) \nonumber \\
&\leq& w(x^0) \left[ 2 C_{n,\alpha} c_0 \int_{\Sigma} \frac{1}{|x^0-y^{\lambda}|^{n+\alpha}} d y + c(x^0) \right]
\label{2.8}
\end{eqnarray}

Since $c(x)$ is bounded from below, to derive a contradiction, it is suffice to show that the integral in the above bracket can be arbitrarily large as $\delta$ becomes sufficiently small. To see this,
let
$$D=\{y|\delta<y_1- x^0_1<1,\; |y'-(x^0)'|<1\},$$
$$s=y_1- x^0_1, \;\; \tau=|y'-(x^0)'|,$$ and $\omega_{n-2}$ be the area of $(n-2)$-dimensional unit sphere. Here we write $x = (x_1, x')$.

 Then we have
\begin{eqnarray}\nonumber
\int_{\Sigma} \frac{1}{|x^0-\tilde{y}|^{n+\alpha}}dy&\geq&\int_D \frac{1}{|x^0-y|^{n+\alpha}}dy\\\nonumber
&=& \int_\delta^1\int_0^1\frac{\omega_{n-2}\tau^{n-2} d\tau}{(s^2+
\tau^2)^{\frac{n+\alpha}{2}}}ds\\ \nonumber
&=& \int_\delta^1\int_0^{\frac{1}{s}}\frac{{\omega_{n-2}(st)^{n-2}}s dt}{s^{n
+\alpha}(1+t^2)^{\frac{n+\alpha}{2}}}ds\\\nonumber
&=& \int_\delta^1\frac{1}{s^{ 1+\alpha}}\int_0^{\frac{1}{s}}\frac{\omega_{n-2}t^{n-2}
dt}{(1+t^2)^{\frac{n+\alpha}{2}}}ds\\\nonumber
&\geq & \int_\delta^1\frac{1}{s^{1+\alpha}}\int_0^1\frac{\omega_{n-2}t^{n-2} dt}{(1+t^2)^{\frac{n+\alpha}{2}}}ds\\\label{8}
&\geq & C\int_\delta^1\frac{1}{s^{1+\alpha}}ds \geq \frac{c}{\delta^{\alpha}}.
\end{eqnarray}

Combining (\ref{2.8}) with (\ref{8}), we arrive at
\begin{eqnarray*}
& &  F(u_{\lambda}(x^0)) - F(u(x^0)) + c(x^0) w(x^0) \nonumber \\
&\leq& w(x^0) \left[ 2 C_{n,\alpha} c_0 \frac{c}{\delta^{\alpha}} + c(x^0) \right]
\end{eqnarray*}

Notice that $w(x^0) < 0$, so for sufficiently small $\delta$ we have
$$ F(u_{\lambda}(x^0)) - F(u(x^0)) + c(x^0) w(x^0) < 0.$$
This contradicts the equation and hence proves the theorem.

\medskip

As one will see from the proof of the above theorem, the contradiction arguments are conducted at a negative minimum of $w$. Hence when working on an unbounded domain, one needs to rule out the possibility that such minima would ``leak'' to infinity. This can be done when $c(x)$ decays
``faster'' than $1/|x|^{\alpha}$ near infinity.

\begin{thm} ( Decay at Infinity.)

Let $\Omega$ be an unbounded region in $\Sigma$.
Assume $w \in L_{\alpha} \cap C_{loc}^{1,1}(\Omega)$ is lower semi-continuous on $\bar{\Omega}$ and  is a solution of
$$
\left\{\begin{array}{ll}
F(u_{\lambda}(x)) - F(u(x)) +c(x)w(x)\geq 0  &\mbox{ in } \Omega,\\
w(x) \geq 0 & \mbox{ in } \Sigma \backslash\Omega,\\
w(x^{\lambda})=-w(x)  & \mbox{ in } \Sigma,
\end{array}
\right.
$$

with
$$
\underset{|x|\rightarrow \infty}{\underline{\lim}}|x|^{\alpha}c(x)\geq0,
$$
then there exists a constant $R_0>0$ ( depending on $c(x)$, but independent of $w$ ), such that if
$$
w(x^0)=\underset{\Omega}{\min}\;w(x)<0,
$$
 then
 $$
|x^0|\leq R_0.
$$
\label{thm2.4}
\end{thm}

{\bf Proof.}

Otherwise, there exists an $x^0 \in \Omega$, such that
$$w(x^0) = \min_{\Omega} w < 0. $$
Again by the {\em key inequality} (\ref{ki}), we deduce
\begin{eqnarray}
& & F(u_{\lambda}(x^0)) - F(u(x^0)) + c(x^0) w(x^0) \nonumber \\
&\leq& w(x^0) \left[ 2 C_{n,\alpha} c_0 \int_{\Sigma} \frac{1}{|x^0-y^{\lambda}|^{n+\alpha}} d y + c(x^0) \right]
\label{2.9}
\end{eqnarray}

We now estimate the above integral. Let $\Sigma^c = \mathbb{R}^n \setminus \Sigma$.
Choose a point in $\Sigma^c$: $x^1=(3|x^0|+x^0_1, (x^0)')$, then
$B_{|x^0|}(x^1)\subset\Sigma^c$. It follows that
\begin{eqnarray*}
\int_\Sigma \frac{1}{|x^0-y^{\lambda}|^{n+\alpha}}dy &=& \int_{\Sigma^c} \frac{1}{|x^0-y|^{n+\alpha}}dy\\
&\geq&\int_{B_{|x^0|}(x^1)} \frac{1}{|x^0-y|^{n+\alpha}}dy\\
&\geq&\int_{B_{|x^0|}(x^1)} \frac{1}{4^{n+\alpha}|x^0|^{n+\alpha}}dy\\
&=&\frac{\omega_n}{4^{n+\alpha}|x^0|^{\alpha}}.
\end{eqnarray*}

Then from the equation and (\ref{2.9}), we have
\begin{eqnarray*}
0&\leq & F(u_{\lambda}(x^0)) - F(u(x^0)) + c(x^0)w(x^0)\\
&\leq&\left[\frac{2\omega_nC_{n,\alpha} c_0}{4^{n+\alpha}|x^0|^{\alpha}}+
c(x^0)\right]w(x^0).
\end{eqnarray*}
Or equivalently,
$$\frac{2 \omega_nC_{n,\alpha} c_0}{4^{n+\alpha} |x^0|^{\alpha}}+c(x^0)\leq 0.$$
Now if $|x^0|$ is sufficiently large, this would contradict the decay assumption on $c(x)$. Therefore, we must have $$w(x) \geq 0, \;\; x \in \Omega.$$ This verifies the theorem.

\section{Applications--Symmetry and Non-existence of Solutions}

\subsection{Symmetry of Solutions In a Unit Ball}

Consider
\begin{equation}
\left\{ \begin{array}{ll}
F_{\alpha}(u(x))  = f(u(x)) , & x \in B_1(0), \\
u(x) = 0 , & x \not{\in} B_1(0).
\end{array}
\right.
\label{fu1}
\end{equation}
We prove
\begin{thm}
Assume that $u \in L_{\alpha}\cap C_{loc}^{1,1}(B_1(0))$ is a positive solution of (\ref{fu1}) with $f(\cdot)$ being Lipschitz continuous. Then $u$ must be radially symmetric and monotone decreasing about the origin.
\label{thmfu}
\end{thm}

{\bf Proof.}  Let $T_{\lambda}, x^{\lambda}, u_{\lambda},$ and $w_{\lambda}$ as defined in the previous section. Let
$$ \Sigma_{\lambda} = \{ x \in B_1(0) \mid x_1 < \lambda \} .$$
Then it is easy to verify that
$$ F_\alpha(u_{\lambda}(x)) - F_\alpha(u(x)) + c_{\lambda}(x) w_{\lambda}(x) = 0 , \;\; x \in \Sigma_{\lambda}, $$
where
$$ c_{\lambda}(x) = \frac{f(u(x)) - f(u_{\lambda}(x))}{u(x) - u_{\lambda}(x)} .$$

Our Lipschitz continuity assumption on $f$ guarantees that $c_{\lambda}(x)$ is uniformly bounded
from below. Now we can apply Theorem \ref{thm2.3} ({\em narrow region principle}) to conclude that
for $\lambda > -1$ and sufficiently close to $-1$,
\begin{equation}
w_{\lambda}(x) \geq 0 , \;\; x \in \Sigma_{\lambda} ;
\label{wgeq01}
\end{equation}
because $\Sigma_{\lambda}$ is a narrow region  for such $\lambda$.

Define
$$ \lambda_0 = \sup \{ \lambda \leq 0 \mid w_{\mu} (x) \geq 0, x \in \Sigma_{\mu}; \mu \leq \lambda \}.$$
Then we must have $$\lambda_0 = 0.$$ Otherwise, suppose that $\lambda_0 < 0$, we show that the plane can be moved to the right a little more and inequality (\ref{wgeq01}) is still valid. More precisely, there exists small $\epsilon >0$, such that for all $\lambda \in [\lambda_0, \lambda_0 + \epsilon)$, inequality (\ref{wgeq01}) holds, which contradicts the definition of $\lambda_0$.

First, since $w_{\lambda_0}(x)$ is not identically zero, from the strong {\em maximum principle} (Theorem \ref{thm2.2}), we have
$$w_{\lambda_0}(x)> 0, \;\; \forall \, x \in \Sigma_{\lambda_0} .$$
Thus for any $\delta >0$,
$$w_{\lambda_0}(x)> c_{\delta} > 0, \;\; \forall \, x \in \Sigma_{\lambda_0 -\delta}.$$

By the continuity of $w_{\lambda}$ with respect to $\lambda$, there exist $\epsilon >0$, such that
\be w_{\lambda}(x) \geq 0, \;\; \forall \, x \in \Sigma_{\lambda_0 -\delta}, \; \forall \lambda \in [\lambda_0, \lambda_0 +\epsilon).
\label{2.5b}
\ee

  In the {\em narrow region principle} (Theorem \ref{thm2.3}), let
$$\Sigma = \Sigma_{\lambda} \mbox{ and the narrow region } \Omega = \Sigma_{\lambda} \setminus \Sigma_{\lambda_0 -\delta},$$
then we have
$$w_{\lambda}(x) \geq 0, \;\; \forall \, x \in \Sigma_{\lambda} \setminus \Sigma_{\lambda_0 -\delta}.$$
This together with (\ref{2.5b}) implies
$$w_{\lambda}(x) \geq 0, \;\; \forall \, x \in \Sigma_{\lambda}, \; \forall \lambda \in [\lambda_0, \lambda_0 +\epsilon).$$
This contradicts the definition of $\lambda_0$. Therefore, we must have $\lambda_0 =0.$ It follows that

$$
w_{0}(x) \geq 0 , \;\; x \in \Sigma_{0} ;
$$
or more apparently,
\begin{equation}
u(-x_1, x_2, \cdots, x_n) \leq u(x_1, x_2, \cdots, x_n) , \;\; 0<x_1<1.
\label{u0}
\end{equation}

Since the $x_1$-direction can be chosen arbitrarily, (\ref{u0}) implies $u$ is radially symmetric about the origin. The monotonicity is a consequence of the fact that
$$w_{\lambda}(x) > 0, \;\; \forall \, x \in \Sigma_{\lambda}$$
 holds for all
$-1<\lambda \leq 0$. This completes the proof of the theorem.

\subsection{Symmetry of Solutions in the Whole Space $\mathbb{R}^n$}

Let $T_{\lambda}, x^{\lambda}, u_{\lambda},$ and $w_{\lambda}$ as defined in the previous section. Let
$$ \Sigma_{\lambda} = \{ x \in \mathbb{R}^n \mid x_1 < \lambda \} .$$

\begin{thm}
Assume that $u \in C^{1,1}_{loc} \cap L_{\alpha}$ is a positive solution of
\be
F_{\alpha}(u(x)) = g(u(x)), \;\; x \in \mathbb{R}^n .
\label{eqws}
\ee

Suppose, for some $\gamma >0$,
\be
u(x)=o(\frac{1}{|x|^{\gamma}}), \mbox{ as } |x|\ra \infty,
\label{u1}
\ee
and
\be
g'(s) \leq s^q , \;\; \mbox{ with } q \gamma \geq \alpha .
\label{g1}
\ee

Then $u$ must be radially symmetric about some point in $\mathbb{R}^n$.
\end{thm}

{\bf Proof.} We carry out the proof in two steps. To begin with, we show that for $\lambda$ sufficiently negative,
\begin{equation}
w_{\lambda} (x) \geq 0 , \;\; \forall x \in \Sigma_{\lambda}
\label{sa21}
\end{equation}
with an application of the \emph{decay at infinity} (Theorem \ref{thm2.4}).

Next, we move the plane $T_{\lambda}$ along the
the $x_{1}$-axis to the right as long as inequality
(\ref{sa21}) holds. The plane will eventually stop at some limiting position at $\lambda=\lambda_o$. Then we claim that
$$ w_{\lambda_o} (x) \equiv 0 , \;\; \forall x \in
\Sigma_{\lambda_o}.$$ The symmetry and monotone decreasing properties of
solution $u$ about $T_{\lambda_o}$ follows naturally from the proof.
Also, because of the arbitrariness of the $x_1$-axis, we conclude that
$u$ must be radially symmetric and monotone about some point.

\emph{Step 1.} \emph{Start moving the plane $T_\lambda$ along the $x_1-$axis
from near $-\infty$ to the right.}

By the {\em mean value theorem} it is easy to see that
\begin{equation}
F_{\alpha}(u_{\lambda}(x)) - F_{\alpha}(u(x)) = g(u_{\lambda}(x)) - g(u(x)) = g'(\psi_{\lambda}(x)) w_{\lambda} (x) .
\label{sa22}
\end{equation}
where $\psi_{\lambda}(x)$ is between $u_{\lambda}(x)$
and $u(x)$. By the \emph{decay at infinity}
argument (Theorem \ref{thm2.4}), it suffices to check the decay
rate of $g'(\psi_{\lambda}(x))$, and to be more precise, only at
the points $\tilde{x}$ where $w_{\lambda}$ is negative. Since
$$u_{\lambda}(\tilde{x}) < u(\tilde{x}),$$
we have
$$0 \leq  u_{\lambda} (\tilde{x}) \leq \psi_{\lambda} (\tilde{x}) \leq u(\tilde{x}) .$$
The decay assumptions (\ref{u1}) and (\ref{g1})
instantly yields that
$$g'(\psi_{\lambda}(\tilde{x})) = o \left(\frac{1}{|\tilde{x}|^\alpha}\right) .$$
Consequently, there exists $R_0 >0$, such that, if $x^o$ is a negative minimum of $w_{\lambda}$ in
$\Sigma_{\lambda}$, then
\be
|x^o| \leq R_0 .
\label{3.5}
\ee

 Now take $\Omega = \Sigma_{\lambda}$, then by Theorem \ref{thm2.4}, it's easy to conclude that for $\lambda \leq - R_0$, we must have (\ref{sa21}). This completes the
preparation for the moving of planes.
\medskip

\emph{Step 2.} \emph{Keep moving the plane to the limiting position
$T_{\lambda_0}$ as long as (\ref{sa21}) holds.}

Let $$\lambda_0=\sup\{\lambda \mid w_\mu(x)\geq0, \; \forall x \in \Sigma_\mu
, \mu\leq\lambda\}.$$
In this part, we will show that
\be
 w_{\lambda_0}(x)\equiv0, \quad x \in \Sigma_{\lambda_0}.\label{sa24}
 \ee
Otherwise, the plane $T_{\lambda_0}$ can still be moved further to the right. More rigorously, there exists a $\delta_0 > 0$
such that for all $ 0 < \delta < \delta_0$, we have
\begin{equation}
w_{\lambda_0 + \delta} (x) \geq 0 \; , \;\; \forall x \in
\Sigma_{\lambda_0 + \delta} . \label{sa25}
\end{equation}
This would contradict the definition of $\lambda_0$, and hence
(\ref{sa24}) must be true. Below we prove (\ref{sa25}).

Suppose (\ref{sa24}) is false, then $w_{\lambda_0}$ is positive somewhere in $\Sigma_{\lambda_0}$, and the {\em strong maximum principle for anti-symmetric functions} (Theorem \ref{thm2.2}) implies
$$w_{\lambda_0}(x) > 0 , \;\; x \in \Sigma_{\lambda_0} .$$

It follows that for any positive number $\sigma$,
\be
w_{\lambda_0}(x)\geq c_o >0, \quad x \in \overline{\Sigma_{\lambda_0-\sigma}
\cap B_{R_0}(0)},
\label{3.6}
\ee
where $R_0$ is defined in {\em Step 1}. Since $w_{\lambda}$ depends on $\lambda$ continuously, for all sufficiently small $\delta>0$, we have
\be
w_{\lambda_0+\delta}(x)\geq0, \quad x \in \overline{\Sigma_{\lambda_0-\sigma}\cap
B_{R_0}(0)}.\label{sa26}
\ee

Suppose (\ref{sa25}) is false, then there exist $x^o \in \Sigma_{\lambda_0 +\delta}$, such that
$$ w_{\lambda_0 + \delta}(x^o) = \min_{\Sigma_{\lambda_0+\delta}} w_{\lambda_0 + \delta} < 0 .$$
From the {\em decay at infinity theorem} and (\ref{sa25}), we must have
\be
 x^o \in
 (\Sigma_{\lambda_0+\delta}\backslash\Sigma_{\lambda_0-\sigma})
\cap B_{R_0}(0).
\label{3.7}
\ee

Notice that $(\Sigma_{\lambda_0+\delta}\backslash\Sigma_{\lambda_0-\sigma})
\cap B_{R_0}(0)$ is a narrow region for sufficiently small $\sigma$ and $\delta$, and by the {\em narrow region principle} (Theorem \ref{thm2.4}),
$w_{\lambda_0 + \delta}$ cannot attain its negative minimum here, which contradicts (\ref{3.7}),
and hence (\ref{sa25}) holds.

 This completes the proof of the Theorem.

\subsection{Non-existence of Solutions on a Half Space}

We investigate a Dirichlet problem on an upper half space
$$\mathbb{R}_+^n = \{ x = (x_1, \cdots, x_n) \mid x_n > 0 \}.$$

Consider
\begin{equation}
\left\{\begin{array}{ll}
F_{\alpha}(u(x))= h(u(x)), & x \in \mathbb{R}_+^n, \\
u(x) \equiv 0 , & x \not{\in} \mathbb{R}_+^n .
\end{array}
\right.
\label{e1h1}
\end{equation}

\begin{thm}
Assume that $u \in L_{\alpha} \cap C^{1,1}_{loc}$ is a nonnegative solution of
problem (\ref{e1h1}). Suppose
\be \lim_{|x| \ra \infty} u(x) = 0,
\label{u2}
\ee
\be
h(s) \mbox{ is Lipschitz continuous in the range of } u, \mbox{ and } h(0)=0.
\label{h1}
\ee

Then $u \equiv 0$.
\label{thmh}
\end{thm}

Based on the assumption (\ref{u2}) and $h(0)=0$, from the proof of Theorem \ref{thm2.1}, one can see that we have
 $$\mbox{ either } u(x) > 0 \mbox{ or } u(x) \equiv 0 , \;\; \forall \, x \in \mathbb{R}^n_+ .$$
 Hence in the following, we may assume that $u >0$ in $\mathbb{R}^n_+.$

 Now we carry on the method of moving planes on the solution $u$ along $x_n$ direction.

 Let $$T_{\lambda} = \{ x \in \mathbb{R}^n \mid x_n = \lambda \} , \;\; \lambda > 0 ,$$
 and $$ \Sigma_{\lambda} = \{ x \in \mathbb{R}^n \mid 0< x_n < \lambda \} . $$
 Let $$ x^{\lambda} = (x_1, \cdots, x_{n-1}, 2\lambda - x_n)$$ be the reflection of $x$ about the plane $T_{\lambda}$. Denote $w_{\lambda} (x) = u(x^{\lambda}) - u(x)$.

 The key ingredient in this proof is the {\em narrow region principle} (Theorem \ref{thm2.3}). To see that it still applies in this situation, we only need to take
 $$\Sigma = \Sigma_{\lambda} \cup \mathbb{R}^n_- ,$$
 where $$ \mathbb{R}^n_- = \{ x \in \mathbb{R}^n \mid x_n \leq 0 \}. $$

 {\em Step 1.} For $\lambda$ sufficiently small, since $\Sigma_{\lambda}$ is a narrow region, we have immediately
 \be w_{\lambda}(x) \geq 0, \;\; \forall \, x \in \Sigma_{\lambda}.
 \label{w1}
 \ee

 {\em Step 2.}  (\ref{w1}) provides a starting point, from which we can move the plane $T_{\lambda}$ upward as long as inequality (\ref{w1}) holds. Define
 $$ \lambda_o = \sup \{ \lambda \mid w_{\mu} (x) \geq 0, x \in \Sigma_{\mu}; \mu \leq \lambda \}.$$
 We show that
 \be \lambda_o = \infty .
 \label{y1}
 \ee

 Otherwise, if $\lambda_o < \infty$, then by (\ref{w1}), combining the {\em narrow region principle} and {\em decay at infinity} and going through the similar arguments as in the previous subsection, we are able to show that
 $$ w_{\lambda_o}(x) \equiv 0 \;\; \mbox{ in } \Sigma_{\lambda_o} ,$$
 which implies
 $$ u(x_1, \cdots, x_{n-1}, 2\lambda_o) = u(x_1, \cdots, x_{n-1}, 0) = 0 . $$
 This is impossible, because we assume that $u>0$ in $\mathbb{R}^n_+.$

 Therefore, (\ref{y1}) must be valid. Consequently, the solution $u(x)$ is monotone increasing with respect to $x_n$. This contradicts (\ref{u1}).

\section{The limit of $F_{\alpha}(u)$ as $\alpha \ra 2$.}

In this section, we let $\alpha \ra 2$ and investigate the limit  of $F_{\alpha}(u(x))$ for each fixed $x$.
\begin{thm}
Assume that $u \in C^{1,1}_{loc} \cap L_{\alpha}$, $G(\cdot)$ is second order differentiable, and $G(0) = 0$.
Then
\be \lim_{\alpha \ra 2} F_{\alpha}(u(x)) = a (-\lap u)(x) + b |\grad u(x)|^2 ,
\label{limit}
\ee
where $a$ and $b$ are constant multiple of $G'(0)$ and $G''(0)$ respectively.
\label{thm4.1}
\end{thm}

{\bf Proof.}
We use the fact $C_{n,\alpha} = c_n (2-\alpha)$ with some constant $c_n$ depending on $n$.

It follows that
\begin{eqnarray}
&& F_{\alpha}(u(x)) \nonumber \\
&=& c_n (2-\alpha) PV \int_{B_{\epsilon}(x)} \frac{G(u(x)-u(y))}{|x-y|^{n+\alpha}} d y +
c_n (2-\alpha)  \int_{\mathbb{R}^n \setminus B_{\epsilon}(x)} \frac{G(u(x)-u(y))}{|x-y|^{n+\alpha}} d y
\nonumber \\
&=& I_1 + I_2 .
\label{4.1}
\end{eqnarray}

First fix $\epsilon$ and let $\alpha \ra 2$. Then obviously
\be
I_2 \ra 0.
\label{4.2}
\ee

To estimate $I_1$, we apply Taylor expansion on $G$ near $0$ and on $u$ near $x$.
\begin{eqnarray*}
&&G(u(x)-u(y))  \\
&=&G'(0) (u(x)-u(y)) + \frac{1}{2} G''(0) (u(x)-u(y))^2 + o(\epsilon) (u(x)-u(y))^2  \\
&=& - G'(0) [ \grad u(x)\cdot z + u_{ij}(x) z_i z_j ] + \frac{1}{2} G''(0) ( u_i(x) u_j(x) z_i z_j )+ o(\epsilon) |z|^2.
\end{eqnarray*}
Here we wrote $z=y-x$ to avoid length expressions, adapted the summation convention that
$$ u_{ij}z_i z_j = \sum_{i,j} u_{ij}z_i z_j ,$$
and $$o(\epsilon) \ra 0, \;\; \mbox{ as } \epsilon \ra 0.$$

Evaluate $I_1$ in four separate parts:
\be
I_1 = I_{11} + I_{12} + I_{13} + I_{14}.
\label{4.4}
\ee

Due to symmetry,
\be
I_{11} = - c_n (2-\alpha) G'(0) PV \int_{B_{\epsilon}(0)} \frac{\grad u(x) \cdot z}{|z|^{n+\alpha}} d z =0.
\label{4.5}
\ee

\begin{eqnarray}
I_{12} &=& -c_n (2-\alpha) G'(0) u_{ij}(x) \int_{B_{\epsilon}(0)} \frac{z_i z_j}{|z|^{n+\alpha}} d z \nonumber \\
&=& -c_n (2-\alpha) G'(0) u_{ii}(x) \int_{B_{\epsilon}(0)} \frac{z_i^2}{|z|^{n+\alpha}} d z \nonumber \\
&=& -c_n (2-\alpha) G'(0) \frac{\lap u(x)}{n} \int_{B_{\epsilon}(0)} \frac{1}{|z|^{n+\alpha-2}} d z \nonumber \\
&=& a (-\lap u(x)) \epsilon^{2-\alpha} \nonumber \\
&\ra& a (-\lap u(x)), \mbox{ as } \alpha \ra 2 .
\label{4.6}
\end{eqnarray}

\begin{eqnarray}
I_{13} &=& 1/2 c_n (2-\alpha) G''(0) u_i(x) u_j(x) \int_{B_{\epsilon}(0)} \frac{z_i z_j}{|z|^{n+\alpha}} d z \nonumber \\
&=& 1/2 c_n (2-\alpha) G''(0) u_i^2(x) \int_{B_{\epsilon}(0)} \frac{z_i^2}{|z|^{n+\alpha}} d z \nonumber \\
&=& 1/2 c_n (2-\alpha) G''(0) \frac{|\grad u(x)|^2}{n} \int_{B_{\epsilon}(0)} \frac{1}{|z|^{n+\alpha-2}} d z \nonumber \\
&=& b (|\grad u(x)|^2) \epsilon^{2-\alpha} \nonumber \\
&\ra& b (|\grad u(x)|^2), \mbox{ as } \alpha \ra 2 .
\label{4.7}
\end{eqnarray}

Then we let $\epsilon \ra 0$.
\be
I_{14} = o(\epsilon) \ra 0 .
\label{4.8}
\ee

Combining (\ref{4.1}), (\ref{4.2}, (\ref{4.4}), (\ref{4.5}), (\ref{4.6}), (\ref{4.7}), and (\ref{4.8}), we prove the theorem.

\bigskip

{\em Authors' Addresses and E-mails:}
\medskip

Wenxiong Chen

Department of Mathematical Sciences

Yeshiva University

New York, NY, 10033 USA

wchen@yu.edu
\medskip

Congming Li

Department of Applied Mathematics

University of Colorado,

Boulder CO USA

cli@clorado.edu
\medskip

Guanfeng Li

Department of Mathematics

Harbin Institute of Technology

Harbin, China

liguanfeng@hit.edu.cn

\end{document}